\documentclass[12pt]{amsart}
\usepackage{color}
\usepackage{epsfig}
\usepackage{verbatim}
\usepackage{clrscode3e}
\usepackage{parskip}

\newtheorem{theorem}{Theorem}[section]
\newtheorem{lemma}[theorem]{Lemma}

\begin{document}

\title {Deletable edges in 3-connected graphs and their applications}
\maketitle

\begin{center}
S. R. Kingan\footnote{Department of Mathematics,
Brooklyn College CUNY, 2900 Bedford Avenue, Brooklyn, NY 11210,
and The Graduate Center, 365 Fifth Avenue, New
York, NY 10016. Email: skingan@brooklyn.cuny.edu.} 
\end{center}
\begin{abstract} Let $G$ and $H$ be simple 3-connected graphs such that $G$ has an $H$-minor. An edge $e$ in $G$ is called {\it $H$-deletable} if $G\backslash e$ is 3-connected and has an $H$-minor. The main result in this paper establishes that, if $G$ has no $H$-deletable edges, then there exists a sequence of simple 3-connected graphs $G_0, \dots , G_k$ with no $H$-deletable edges such that $G_0\cong H$, $G_k= G$, and for $1 \le i \le k$ one of three possibilities holds: $G_{i-1}= G_i/f$; $G_{i-1}=G_i/f \backslash e$ where $e$ and $f$ are incident to a degree 3 vertex in $G_i$; or $G_{i-1}=G_i-w$ where $w$ is a degree $3$ vertex in $G_i$.
Several applications are given including a graph theoretic proof of the matroid theory result known as the Strong Splitter Theorem, a short new proof of Dirac's characterization of 3-connected graphs with no minor isomorphic to the prism graph, and an extension of a result by Halin that bounds the number of edges in a minimally 3-connected graph. Halin proved that if $G$ is a minimally $3$-connected graph on $n\ge 8$ vertices, then $|E(G)|\le 3n-9$ and equality holds if and only if $G\cong K_{3, n-3}$. We give a different proof of Halin's result and extend it by identifying the minimally 3-connected infinite family of graphs with $|E(G)|=3n-10$. \end{abstract} 

\section{\bf Introduction}

Suppose $G$ and $H$ are simple 3-connected graphs such that $H$ is a proper minor of $G$. In 1980 Seymour gave a remarkable theorem for reducing $G$ to $H$ using deletions and contractions while preserving 3-connectivity and the $H$-minor at each step \cite{Seymour1980}. He proved that, except for a few types of graphs, there is an edge $e$ such that either $G\backslash e$ or $G/e$ is simple 3-connected and has an $H$-minor. An edge $e$ in $G$ is called {\it $H$-deletable} if $G\backslash e$ is 3-connected and has a minor isomorphic to $H$. In this paper we will prove that if $G$ has no $H$-deletable edges, then there exists a simple $3$-connected graph $G'$ with an $H$-minor, but no $H$-deletable edges, such that one of three possibilities holds: $G'= G/f$; $G'= G/f \backslash e$ where edges $e$ and $f$ are incident to a degree $3$ vertex in $G$; or $G'=G-w$ where $w$ is a degree $3$ vertex in $G$. 
As a consequence, we get a sequence of simple 3-connected graphs $G_0, \dots , G_k$ with an $H$-minor, but no $H$-deletable edges, such that $G_0\cong H$, $G_k= G$, and for $1 \le i \le k$ one of the previously mentioned three possibilities holds. When such a result is reversed it allows for the systematic construction of an isomorphic copy of $G$ from $H$ and has applications to excluded minor results. Furthermore, if $G$ has $H$-deletable edges, then the number of $H$-deletable edges is $$|E(G)|-|E(H)|-3(|V(G)|-|V(H)|).$$

An edge $e$ in a 3-connected graph $G$ is called {\it contractible} if $G/ e$ is 3-connected (and not necessarily simple). Considerable research has been done on the number and distribution of contractible edges beginning with Tutte's 1961 result that a 3-connected graph has a contractible edge \cite{Tutte1961}. For example, it is proved in \cite{AndoEnomotoSaito1987} that a 3-connected graph with $n\ge 5$ vertices has at least $\lceil \frac{n}{2}\rceil$ contractible edges. Let $G$ and $H$ be simple 3-connected graphs where $H$ is a proper minor of $G$. It is proved in \cite{Costalonga2020} that $G$ has a forest $F$ with at least $\frac{1}{2}(|E(G)|-|E(H)|)$ edges such that for every $e\in E(F)$, $G/e$ is 3-connected and has an $H$-minor. See \cite{Kawarabayashi2002}, \cite{Kriesell2002}, \cite{Ota1988}, \cite{ReidWu2000}, \cite{Saito1990}, and the references in these papers for additional results. An edge $e$ of a 3‐connected graph $G$ is said to be {\it removable} if $G\backslash e$ is a subdivision of a 3‐connected graph. In \cite{Su1999}, Su proved that every 3-connected graph $G$ with $|V(G)|\ge 5$, except $W_5$ and $W_6$, has at least $\frac{1}{7}(3|V(G)|+18)$ removable edges. See also \cite{Holtonetal1990} and \cite{XuandGuo2019} for results on removable edges. 

This paper introduces $H$-deletable edges in graphs and highlights their usefulness for excluded minor results. We will use deletable edges to give a graph theoretic proof of the matroid result known as the Strong Splitter Theorem, a new proof of an excluded minor result by Dirac, and extend a result of Halin on minimally 3-connected graphs. 

\begin{figure}[h]
\centering
\includegraphics[width=6in]{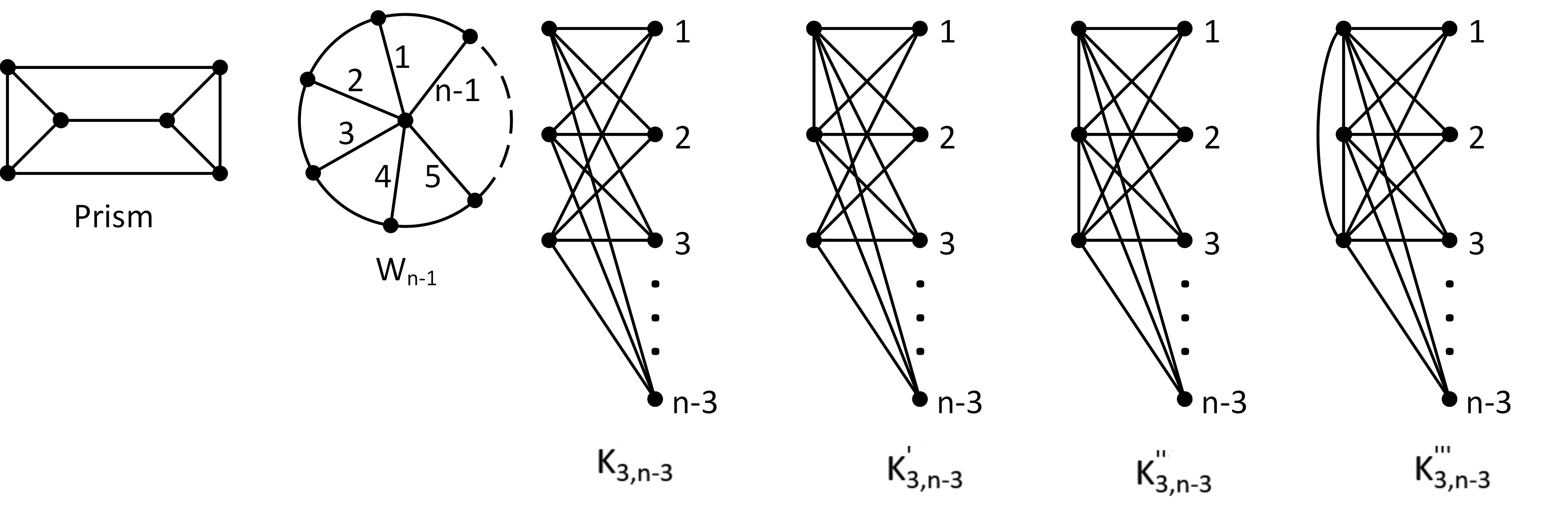}
\caption{Prism, $W_{n-1}$ for $n\ge 4$, and $K_{3, n-3}$, $K'_{3, n-3}$, $K''_{3, n-3}$, $K'''_{3, n-3}$ for $n\ge 6$}
\label{fig-1-EXprism} 
\end{figure}

The wheel graph with $n$ vertices and $n-1$ spokes is denoted by $W_{n-1}$ where $n\ge 4$. The prism graph is the geometric dual of the complete graph $K_5$ with an edge deleted. In 1963, Dirac proved that a simple $3$-connected graph $G$ has no prism minor if and only if $G$ is $K_5$, $K_5\backslash e$, $W_{n-1}$ where $n\ge 4$, or $K_{3,n-3}$, $K'_{3,n-3} $, $K''_{3,n-3}$, or $K'''_{3,n-3}$ where $n\ge 6$ \cite{Dirac1963}. The prism graph and the infinite families are shown in Figure \ref{fig-1-EXprism}. 

\begin{figure}[h]
\centering
\includegraphics[width=0.75in]{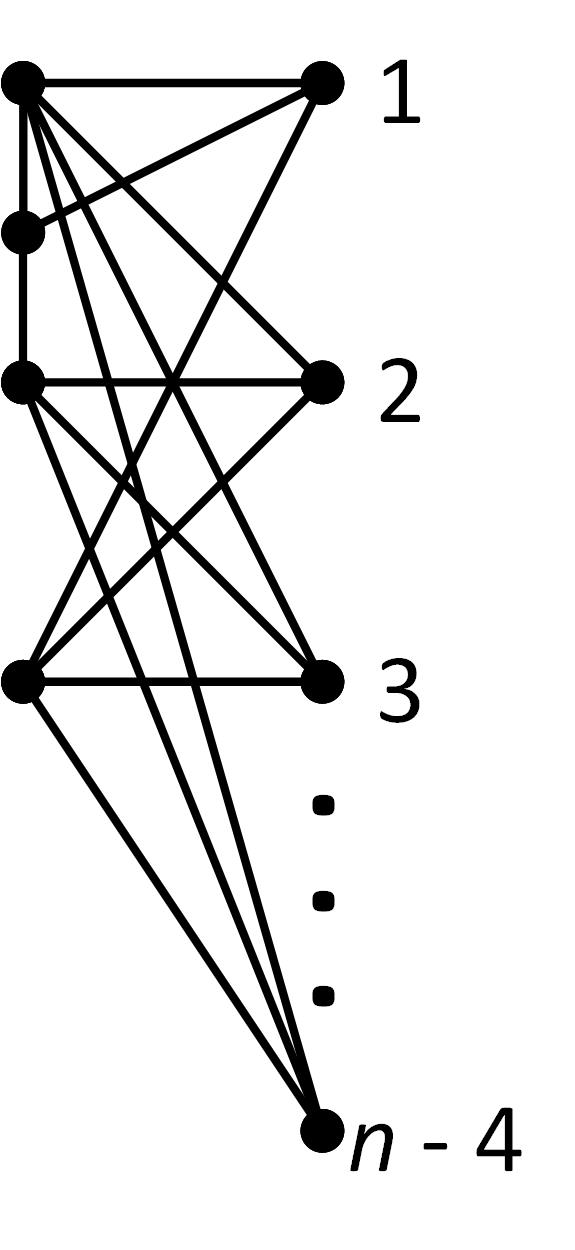}
\caption{The size $3n-10$ minimally 3-connected graph $K^{\perp}_{3,n-4}$}
\label{fig-2-newfamily}
\end{figure}

Dirac's original proof is rather long and complicated \cite{Dirac1963}. In \cite{Brown1965} Brown gave a proof using Tutte's Wheels Theorem \cite{Tutte1961}. Note that these proofs are expressed in terms of finding the 3-connected graphs without two vertex disjoint cycles. 
While updating and further shortening the proof of Dirac's result, we discovered an extension to a result by Halin. In 1969, Halin proved that if $G$ is a minimally $3$-connected graph on $n\ge 8$ vertices, then $|E(G)|\le 3n-9$. Equality holds if and only if $G$ is $K_{3, n-3}$ \cite{Halin1969a}, \cite{Halin1969b}. Halin's result was extended by Mader to $k$-connected graphs \cite{Mader1971}, so this result is often credited to Mader. We will prove that if $G$ is minimally 3-connected with $|E(G)|=3n-10$, then $G$ is the infinite family $K^{\perp}_{3, n-4}$ shown in Figure \ref{fig-2-newfamily}.

Section 2 has terminology and previous results, Section 3 has a construction theorem for graphs without deletable edges, and Section 4 has all the applications.


\section{\bf Terminology and prior results}

\begin{figure}[h]
\centering
\includegraphics[width=3.25in]{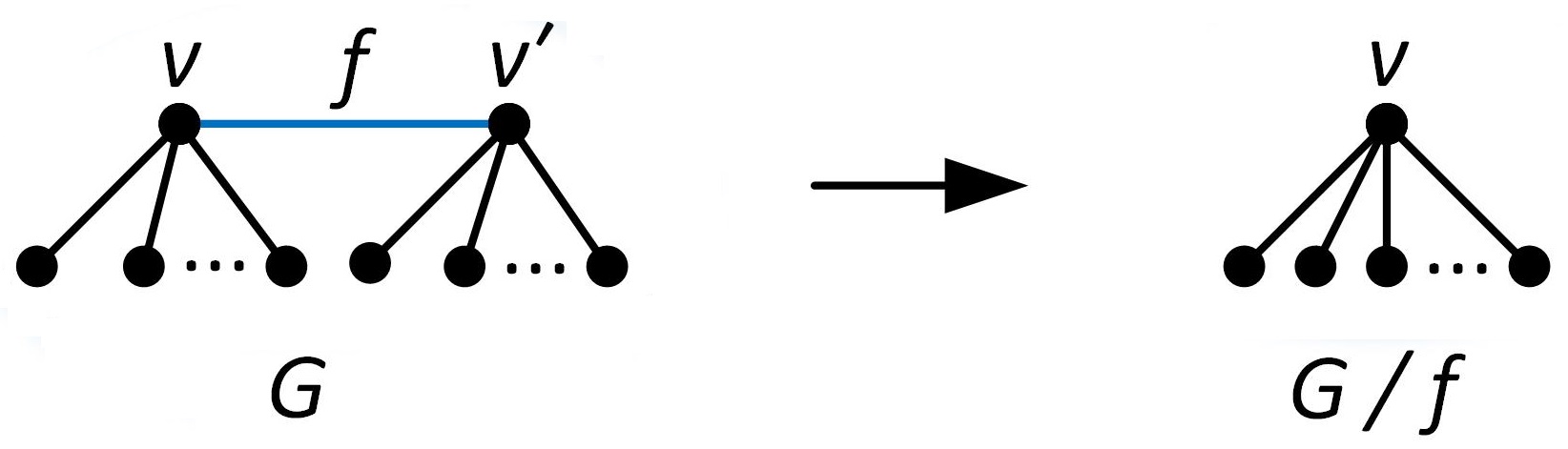}
\caption{Contracting an edge \label{edgecontraction}}
\end{figure}

The terminology and notation follow \cite{Kingan2022} for the most part, except that a graph may have loops and multiple edges. If it does not, then it is called {\it simple}. A graph with at least 4 vertices is {\it $3$-connected} if at least 3 vertices must be removed to disconnect it. By convention $K_4$ is considered 3-connected. To delete an edge $e$ remove it from the graph leaving its end vertices intact. The resulting graph, denoted by $G\backslash e$, is called an {\it edge-deletion} of $G$. 
An edge $e$ in a 3-connected graph $G$ is called {\it deletable} if $G\backslash e$ is 3-connected. A 3-connected graph is {\it minimally $3$-connected} if it has no deletable edges. As noted in the introduction, an edge $e$ in $G$ is $H$-deletable if $G\backslash e$ is 3-connected and has an $H$-minor. 
To contract an edge $f$ with end vertices $v$ and $v'$, collapse the edge by identifying $v$ and $v'$ as one vertex, and delete the loop formed. The resulting graph, denoted by $G / f$, is called an {\it edge-contraction}. Note that if $G$ is simple, then $G\backslash e$ remains simple, but $G/e$ is not necessarily simple. A graph $H$ is a {\it minor} of a graph $G$ if $H$ can be obtained from $G$ by deleting or contracting edges. Figure \ref{edgecontraction} viewed from left to right illustrates the edge-contraction operation. 
Tutte proved that the wheel graph is the unique 3-connected graph such that for every edge $e$, deleting $e$ gives a graph that is not 3-connected and contracting $e$ gives a graph that is not simple and 3-connected \cite{Tutte1961}. Therefore, if $G$ is a simple 3-connected graph that is not a wheel, then there is an edge such that either $G\backslash e$ or $G/e$ is simple and 3-connected. Equivalently, every 3-connected graph, except for wheels, can be constructed from wheels by adding edges and splitting vertices. This is known as the Wheels Theorem. Seymour strengthened Tutte's result by showing that wheels may be replaced by any simple 3-connected proper minor $H$ of $G$, provided that if $H$ is a wheel, then $G$ has no larger wheel minor \cite{Seymour1980}. Seymour's result was in terms of matroids. This is known as the Splitter Theorem. Two years later Negami gave a graph theoretic proof \cite{Negami1982}. Subsequently, Oxley and Coullard simplified the hypothesis to only require that $G$ is not a wheel and $H$ is not $W_3$ \cite{CoullardOxley1992}. See also \cite{Oxley2012} (Corollary 12.3.1). 
\begin{theorem} \label{SplitterTheorem} Suppose $G$ and $H$ are simple $3$-connected graphs such that 
$G$ has a proper $H$-minor, $G\not\cong W_{n-1}$, and $H\not\cong W_3$. Then there exists an edge $e$ such that either $G\backslash e$ or $G/e$ is simple, $3$-connected, and has an $H$-minor. \end{theorem} 

The Wheels Theorem and the Splitter Theorem can also be presented in the bottom-up or constructive format. A graph $G$ with an edge added between non-adjacent vertices is denoted by $G+e$ and called an {\it edge addition} of $G$. Adding an edge between non-adjacent vertices gives exactly one graph. To {\it split} a vertex $v$ with $deg_G(v) \ge 4$, first replace $v$ with two vertices $v_1$ and $v_2$ and a new edge $f=v_1v_2$. Then assign the edges originally incident to $v$ to $v_1$ and $v_2$ so that both $v_1$ and $v_2$ have degree at least 3. The resulting graph is called a {\it vertex split} of $G$ and is denoted by $G\circ f$. Unlike an edge addition, which gives precisely one graph when an edge is added to a pair of non-adjacent vertices, splitting a vertex can give several non-isomorphic graphs depending on the assignment of edges incident to $v_1$ and $v_2$. Observe that if $G$ is 3-connected, then $G+e$ and $G\circ f$ are also 3-connected. 
The constructive format of Theorem \ref{SplitterTheorem} states that an isomorphic copy of $G$ can be constructed from $H$ by a finite sequence of edge additions or vertex splits.
Thus from the top-down perspective, Theorem \ref{SplitterTheorem} implies that there exists a sequence of simple $3$-connected graphs $G_0, \dots ,G_k$ such that $G_0 \cong H$, $G_k=G$, and for $1\le i \le k$, $G_{i-1}= G_i\backslash e$ or $G_{i-1}= G_i/e$, for some edge $e$ in $G_i$. From the constructive perspective, there exists a sequence of simple $3$-connected graphs $G_0, \dots ,G_k$ such that $G_0 = H$, $G_k\cong G$, and for $1\le i \le k$, $G_i= G_{i-1}+e$ or $G_i= G_{i-1}\circ e$

\section{\bf Graphs without $H$-deletable edges}
This section has the statement and proof of a splitter theorem for graphs without $H$-deletable edges.
\begin{theorem} \label{FCT-TopDown} Suppose $G$ and $H$ are simple $3$-connected graphs such that $G$ has a proper $H$-minor, but no $H$-deletable edges. Then there exists a simple $3$-connected graph $G'$ with an $H$-minor, but no $H$-deletable edges, such that: 
\begin{itemize}
\item [(a)] $G'= G/f$; 
\item [(b)] $G'= G/f \backslash e$, where edges $e$ and $f$ are incident to a degree $3$ vertex in $G$; or 
\item [(c)] $G'=G-w$, where $w$ is a degree $3$ vertex in $G$.
\end{itemize}
Moreover, $|V(G')|= |V(G)|-1$ and $|E(G')|\ge |E(G)|-3$.
\end{theorem}

\begin{proof} Suppose $G \cong W_{n-1}$. The only 3-connected minor of a wheel on $n$ vertices is a wheel on fewer vertices and $W_{n-1}/r\backslash s \cong W_{n-2}$, where $r$ is a rim edge and $s$ is a spoke edge. Thus $G$ satisfies (b) in the statement of the theorem. Moreover, $|V(W_{n-2})|= |V(W_{n-1})|-1$ and $|E(W_{n-2})|=|E(W_{n-1})|-2$.
For the rest of the proof, we may assume that $G \not \cong W_{n-1}$ and $H\not\cong W_3$. Since $G$ has an $H$-minor, but no $H$-deletable edges, Theorem \ref{SplitterTheorem} implies that there exists an edge $f$ such that $G/f$ is simple, 3-connected, and has an $H$-minor. Observe that $f$ is not in any triangle of $G$. Now $G/f$ may have $H$-deletable edges. Remove as many edges as needed to obtain a minimally 3-connected graph. Let $X$ be a maximal set of $H$-deletable edges. Then $G'=G/f\backslash X$ is simple, 3-connected, and has an $H$-minor, but no $H$-deletable edges. If $|X|=0$, then we have (a) in the statement of the theorem. Therefore suppose $|X|\ge 1$.

Let $v$ and $v'$ be the end vertices of $f$ in $G$ that are contracted to become a single vertex $v$ in $G/f$. First, we will prove that in $G/f$, every edge in $X$ is incident to $v$. Suppose this is not true for some edge $a\in X$. Since $G/f$ is simple and 3-connected, both end vertices of $a$ have degree at least 4 in $G/f$. When the contraction and deletion operations that form $G/f\backslash X$ are reversed to obtain $G$, this edge $a$ and both its end vertices will remain undisturbed, and therefore will have degree at least 4 in $G$. Since $a$ is an $H$-deletable edge in $G/f$, it will remain an $H$-deletable edge in $G$; a contradiction since $G$ has no $H$-deletable edges.

Next, consider an edge $e$ in $X$. By the argument in the previous paragraph, one end vertex of $e$ in $G/f$ is $v$. Let the other end vertex be $w$. Since $G/f$ is simple and 3-connected, and $e$ is an $H$-deletable edge, $deg_{G/f}(v)\ge 4$ and $deg_{G/f}(w)\ge 4$. In $G$, vertex $w$ is unaffected by the splitting operation and will continue to have degree at least 4. Furthermore, in $G$, at least two of the edges incident to $v$ are unlinked from $v$ and linked to $v'$. Without loss of generality, let $e$ be unlinked from $v$ and linked to $v'$. Observe that $deg_G(w)\ge 4$ and $e$ is a deletable edge in $G/f$. If $deg_G(v')$ is also at least 4, then $e$ remains a deletable edge in $G$; a contradiction. Therefore, $deg_G(v')=3$. Since $f$ is already incident to $v$ and $v'$, we may conclude that $|X|\le 2$.

If $X=\{e\}$, then $G/f\backslash e$ is simple, 3-connected, and has an $H$-minor, but no $H$-deletable edge. Observe that $e$ and $f$ are incident to $v'$ in $G$ and $deg_G(v')=3$. This is the situation in (a). 

If $X=\{e_1, e_2\}$, then $G/f\backslash \{e_1, e_2\}$ is simple, 3-connected, and has an $H$-minor, but no $H$-deletable edge. Both edges $e_1$ and $e_2$ are incident to $v$ in $G/f$ and subsequently to $v'$ in $G$. Observe that all three new edges, $e_1$, $e_2$, and $f$, are incident to the degree 3 vertex $v'$ in $G$ and $G/f\backslash \{e_1, e_2\}=G-v'$. In other words, $G'$ is obtained from $G$ by removing a degree 3 vertex which is (c) in the statement of the theorem. 

Finally, since $G=G/f\backslash X$ and $|X|\le 2$, $|V(G')|= |V(G)|-1$ and $|E(G')|\ge |E(G)|-3$. \end{proof}

\subsection{\bf Constructive format of the main result}
$ $ 

Theorem \ref{FCT-TopDown} may be used to characterize excluded-minor classes, but for this we must look at it from a constructive approach. The set of three edges incident to a degree 3 vertex is called a {\it triad}. 
\begin{theorem} \label{FCT-BottomUp} Suppose $G\not\cong W_{n-1}$ and $H\not\cong W_3$ are simple $3$-connected graphs such that $G$ has a proper $H$-minor, but no $H$-deletable edges. Then there exists a simple $3$-connected graph $G'$ with an $H$-minor, but no $H$-deletable edges, such that: 
\begin{itemize}
\item [(i)] $G = G'\circ f$; 
\item [(ii)] $G = (G'+e) \circ f$, where $e$ and $f$ are in a triad of $G$; or 
\item [(iii)] $G=G'+\{e_1, e_2\}\circ f$, where $\{e_1, e_2, f\}$ is a triad of $G$.
\end{itemize}
Moreover, $|V(G)|=|V(G')+1$ and $|E(G)|\le |E(G')|+3$.
\end{theorem}

\begin{proof} Theorem \ref{SplitterTheorem} implies that we can construct a graph isomorphic to $G$ from $H$ by a sequence of edge additions and vertex splits. 
Since $G$ has no $H$-deletable edges, the last operation in forming $G$ is splitting a vertex. Let $v$ be the vertex in $G'$ that is split to form two vertices $v$ and $v'$ and let $f$ be the edge joining $v$ and $v'$. Then $G= G^+ \circ f$, for some 3-connected graph $G^+$ with $|V(G)|-1$ vertices and an $H$-minor. Now $G^+$ may have $H$-deletable edges. Remove as many edges as needed to obtain a 3-connected graph $G'$ with an $H$-minor and no $H$-deletable edges. Let $X$ be a maximal set of $H$-deletable edges and let $G'=G^+\backslash X$. Then $G= G'+X\circ f.$ As in the proof of Theorem \ref{FCT-TopDown}, $|X|\le 2$. If $|X|=0$, then we have the situation in (i). If $X=\{e\}$, then $G = (G'+e) \circ f$ where $e$ and $f$ are incident to the new degree 3 vertex $v'$. In this case $e$ and $f$ are in a triad of $G$. If $X=\{e_1, e_2\}$, then $G=G'+\{e_1, e_2\}\circ f$. In this case, $e_1$ and $e_2$ are incident with the degree 3 vertex $v'$ in $G$, so $\{e_1, e_2, f\}$ is a triad.
\end{proof}

\begin{figure}[h]
\centering
\includegraphics[width=3.5in]{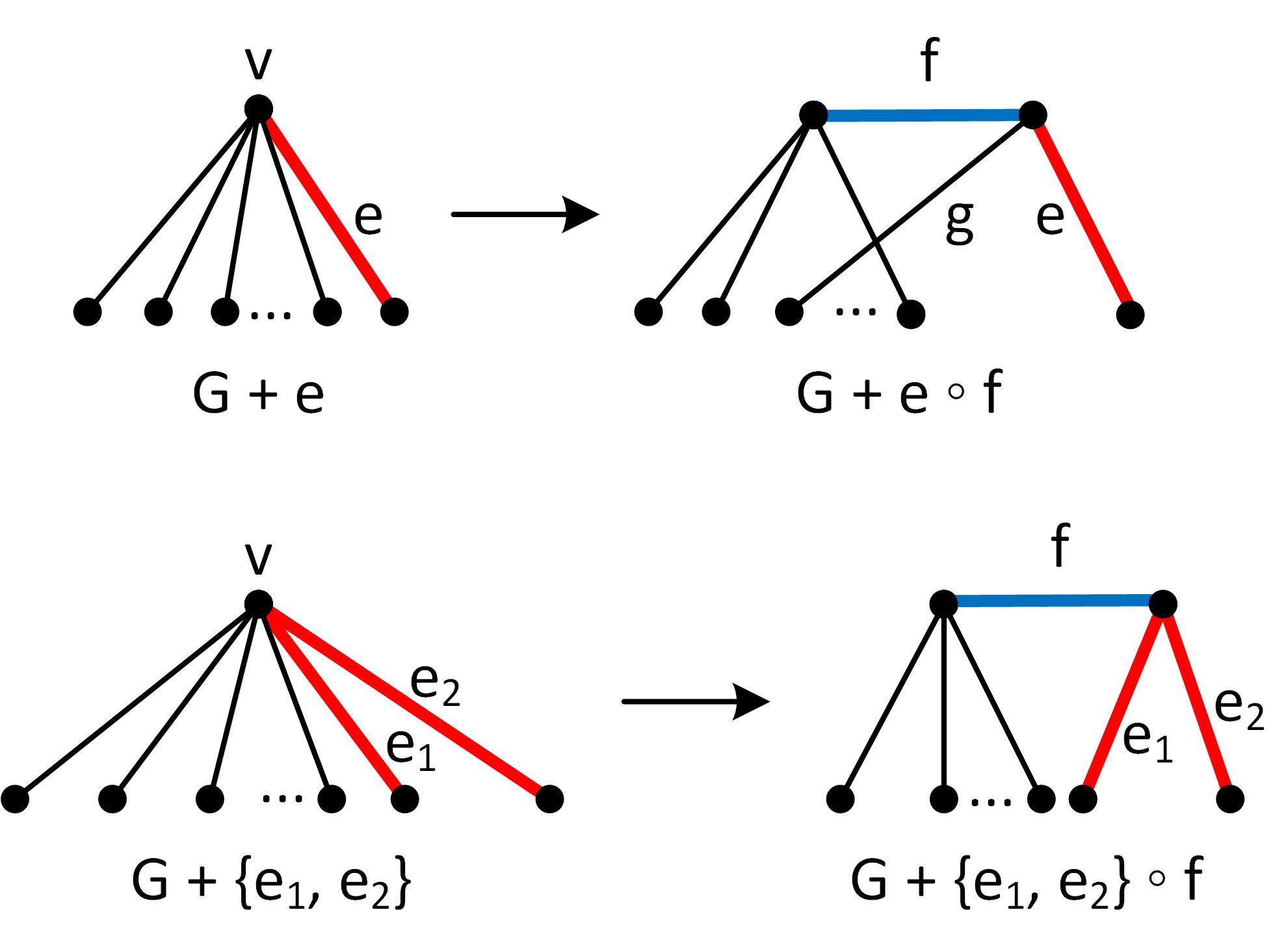}
\caption{Diagram for Theorem \ref{FCT-BottomUp}
\label{FCTFigure}}
\end{figure}

Figure \ref{FCTFigure} illustrates Operations (ii) and (iii) in Theorem \ref{FCT-BottomUp}.

While the top-down version presented in Theorem \ref{FCT-TopDown} and the constructive version presented in Theorem \ref{FCT-BottomUp} are equivalent (except for the exclusion of $W_{n-1}$), practically the latter has efficient graph generation and pattern detection advantages as we will see in subsequent sections. For example, it is not just the case that a degree 3 vertex is added to the graph as stated in Theorem \ref{FCT-TopDown}(iii). Rather, two edges incident to a common vertex $v$ are added to $G$ to form $G+\{e_1, e_2\}$ and the common vertex $v$ is split to form a degree 3 vertex $v'$ in $G$ in precisely one way. Each execution of Theorem \ref{FCT-BottomUp}(iii) gives precisely one graph.


\section{\bf Applications}

\subsection{A splitter theorem for minimally 3-connected graphs}

\begin{theorem} \label{FCT-M3C} Suppose $G$ and $H$ are minimally $3$-connected graphs such that $G$ has a proper $H$-minor. Then there exists a minimally $3$-connected graph $G'$ with an $H$-minor such that: 
\begin{itemize}
\item [(i)] $G'= G/f$; 
\item [(ii)] $G'= G/f \backslash e$, where edges $e$ and $f$ are incident to a degree $3$ vertex in $G$; or 
\item [(iii)] $G'=G-w$, where $w$ is a degree $3$ vertex in $G$.
\end{itemize}
Moreover, $|V(G')|= |V(G)|-1$ and $|E(G')|\ge |E(G)|-3$.
\end{theorem}

\begin{proof} Observe that if $G$ and $H$ are minimally 3-connected graphs such that $G$ has an $H$-minor, then $G$ has no $H$-deletable edges and the result follows from Theorem \ref{FCT-TopDown} \end{proof}

Thus by contracting edges and removing degree 3 vertices, as described in Operations (ii) and (iii), minimal 3-connectivity and the minor is preserved. 


\subsection{Maximum size of a simple 3-connected graph with no H-deletable edge}
$ $

\begin{theorem} \label{FCT-Bound} Suppose $G$ and $H$ are simple $3$-connected graphs such that $G$ has an $H$-minor, but no $H$-deletable edges. Then $|E(G)| \le |E(H)|+ 3(|V(G)|-|V(H)|).$ 
\end{theorem} 

\begin{proof} The proof is by induction on $|V(G)|$. If $|V(G)|=|V(H)|$, then since $G$ has no $H$-deletable edges, $G\cong H$ and the result holds trivially. Suppose $|V(G)|>|V(H)|$. Then Theorem \ref{FCT-BottomUp} implies that $G$ is obtained from a simple $3$-connected graph $G'$ with an $H$-minor, but no $H$-deletable edges, such that $|V(G)|=|V(G')|+1$ vertices and $|E(G)|\le |E(G')|+3$. By the induction hypothesis, 
\begin{align*}
|E(G)| \le |E(G')|+3 &\le |E(H)|+3(|V(G')|-|V(H)|)+3 \\ 
&= |E(H)|+3(|V(G)|-1-|V(H)|)+3 \\
&= |E(H)|+3(|V(G)|-|V(H)|).
\end{align*}
\end{proof}

As a consequence of Theorem \ref{FCT-Bound}, if $|E(G)| > |E(H)|+ 3(|V(G)|-|V(H)|)$, then $G$ has $H$-deletable edges. The number of $H$-deletable edges is $$|E(G)|-|E(H)|- 3(|V(G)|-|V(H)|).$$

\subsection{\label{SST-Subsection} Strong Splitter Theorem for graphs}

$ $

The {\it rank} of a graph $G$, denoted by $r(G)$, is the size of a spanning tree. The Splitter Theorem (Theorem \ref{SplitterTheorem}) implies that we can construct an isomorphic copy of $G$ starting with $H$ and performing a sequence of edge additions and vertex splits. There is no condition on how many edge additions may occur before a vertex split must occur. The Strong Splitter Theorem \cite{KinganLemos2014} implies that at most two consecutive edge additions may be performed in the sequence before a vertex split must be performed, unless the rank of the graphs involved is the same as the rank of $G$. Moreover, when two consecutive edges are added and a vertex split is performed, the three edges form a triad. It is a matroid result and the proof in \cite{KinganLemos2014} is entirely in the language of matroids. Adding an edge to an $n$-vertex graph $G$ intrinsically assumes that there is a larger graph, namely $K_n$, containing $G$ as a deletion-minor from which edges missing in $G$ may be added to $G$ to ``grow" $G$ up to $K_n$. An $n$-element matroid does not necessarily sit inside a larger matroid. As such, the original statement of the Strong Splitter Theorem in \cite{KinganLemos2014} side-steps both the top-down and the constructive approach by talking only about the rank and number of elements of the matroid. 

Suppose $G\not \cong W_{n-1}$ and $H\not \cong W_3$ are simple $3$-connected graphs such that $G$ has an $H$-minor. Further, suppose $t=r(G)-r(H)$. In the language of graphs, the Strong Splitter Theorem states that there exists a sequence of $3$-connected graphs $G_0, G_1,\dots,G_k$ where $k\ge t$, such that:
\begin{enumerate}
\item[(i)] $G_0\cong H$;
\item[(ii)] $G_k=G$;
\item[(iii)] For $1\le i \le t$, $r(G_i)-r(G_{i-1})=1$ and $|E(G_i)|-|E(G_{i-1})|\le 3$; and
\item[(iv)] For $t< i \le k$, $r(G_i)=r(G)$ and $|E(G_i)|-|E(G_{i-1})|=1$.
\end{enumerate}
Moreover, when $|E(G_i)|-|E(G_{i-1})|=3$, then $E(G_i)-E(G_{i-1})$ is a triad of $G_i$.

Combining Theorems \ref{FCT-BottomUp} and \ref{FCT-Bound} gives the following structural result for 3-connected graphs that is equivalent to the Strong Splitter Theorem with more details customized for graphs.

\begin{theorem} \label{SST-BottomUp} Suppose $G\not \cong W_{n-1}$ and $H\not \cong W_3$ are simple $3$-connected graphs such that $G$ has a proper $H$-minor and suppose $G$ has $n$ vertices and $m$ edges where $n\ge |V(H)|+1$ and $m\ge |E(H)|+3$. Then there exists a sequence of simple $3$-connected graphs $G_0, G_1, \dots , G_m$ with an $H$-minor, but no $H$-deletable edges, such that $G_0=H$, $G_n\cong G\backslash D$, $G_m\cong G$, and 
\begin{itemize}
\item [(a)] For $1\le i \le n$:
\begin{itemize}
\item [(i)] $G_i= G_{i-1}\circ f$; 
\item [(ii)] $G_i = (G_{i-1}+e) \circ f$ where $e$ and $f$ are in a triad of $G_i$; or 
\item [(iii)] $G_i=G_{i-1}+\{e_1, e_2\}\circ f$ where $\{e_1, e_2, f\}$ is a triad of $G_i$.
\end{itemize}
\item [(b)] For $n < i \le m$, $G_i= G_{i-1}+e$.
\end{itemize} 
Moreover, $|V(G_i)|=|V(G_{i-1})+1$ and $|E(G_i)|\le |E(G_{i-1})|+3$.
\end{theorem}

\begin{proof} Let $D$ be the set of $H$-deletable edges in $G$. Then $G\backslash D$ is a simple 3-connected graph with an $H$-minor, but no $H$-deletable edges and $V(G\backslash D)= V(G)$. 
Repeated application of Theorem \ref{FCT-BottomUp} implies that there is a sequence of 3-connected graphs $G_0, \dots , G_n$ with an $H$-minor, but no $H$-deletable edges, such that $G_0= H$, $G_n \cong G\backslash D$, and for $1\le i \le n$, $G_i= G_{i-1}\circ f$; $G_i = (G_{i-1}+e) \circ f$ where $e$ and $f$ are in a triad of $G_i$; or $G_i=G_{i-1}+\{e_1, e_2\}\circ f$ where $\{e_1, e_2, f\}$ is a triad of $G_i$. Finally, since $V(G\backslash D)= V(G)$, for 
$n < i \le m$, $G_i$ is obtained from $G_{i-1}$ by adding an edge. 
\end{proof}

\subsection{\label{Dirac} A short new proof of Dirac's Theorem}
$ $

Suppose $J\not \cong W_3$ is a 3-connected graph and let $\mathcal M$ denote the class of graphs with no minor isomorphic to $J$. The goal in an excluded minor result is to determine the graphs in $\mathcal M$ as precisely as possible. It is well known that if $G$ is connected, but not 3-connected, then it can be constructed from its 3-connected proper minors using the operations of 1-sum and 2-sum. Hence, when characterizing excluded minor results, the focus is on finding members in $\mathcal M$ whose connectivity is at least 3. Although Robertson and Seymour developed a structural characterization for $\mathcal M$ in \cite{RobertsonSeymour2003}, this landmark result is an existence result. Practically speaking, few excluded minor results are known \cite{DingLiu2012}. One of the earliest and most well known results is Dirac's characterization of 3-connected graphs with no prism minor \cite{Dirac1963}.

We will call an $n$-vertex 3-connected graph in $\mathcal M$ that has no further edge additions in $\mathcal M$ a {\it monarch} on $n$ vertices. Monarch is a short name for a maximal 3-connected member of $\mathcal M$. We avoid using the word ``extremal'' since we are talking only about 3-connected simple graphs. Before finding the monarchs in $\mathcal M$, we must find the ``minimal'' 3-connected members on $n$ vertices with respect to a 3-connected graph $J \in \mathcal M$. In other words, we must find the 3-connected members of $\mathcal M$ with no $J$-deletable edges. Usually, $J$ is a small 3-connected member in $\mathcal M$. To illustrate this strategy, we give a new short proof of Dirac's characterization of graphs with no prism minor that essentially reduces the entire proof to a ``proof by picture.''

\begin{theorem} \label{DiracTheorem1963} {\bf (Dirac, 1963)} A simple $3$-connected graph $G$ has no prism minor if and only if $G$ is isomorphic to $K_5\backslash e$, $K_5$, $W_{n-1}$, for $n\ge 4$, $K_{3,n-3}$, $K'_{3,n-3} $, $K''_{3,n-3}$, or $K'''_{3,n-3}$, for $n\ge 6$. 
\end{theorem}

\begin{figure}[h]
\centering
\includegraphics[width=6in]{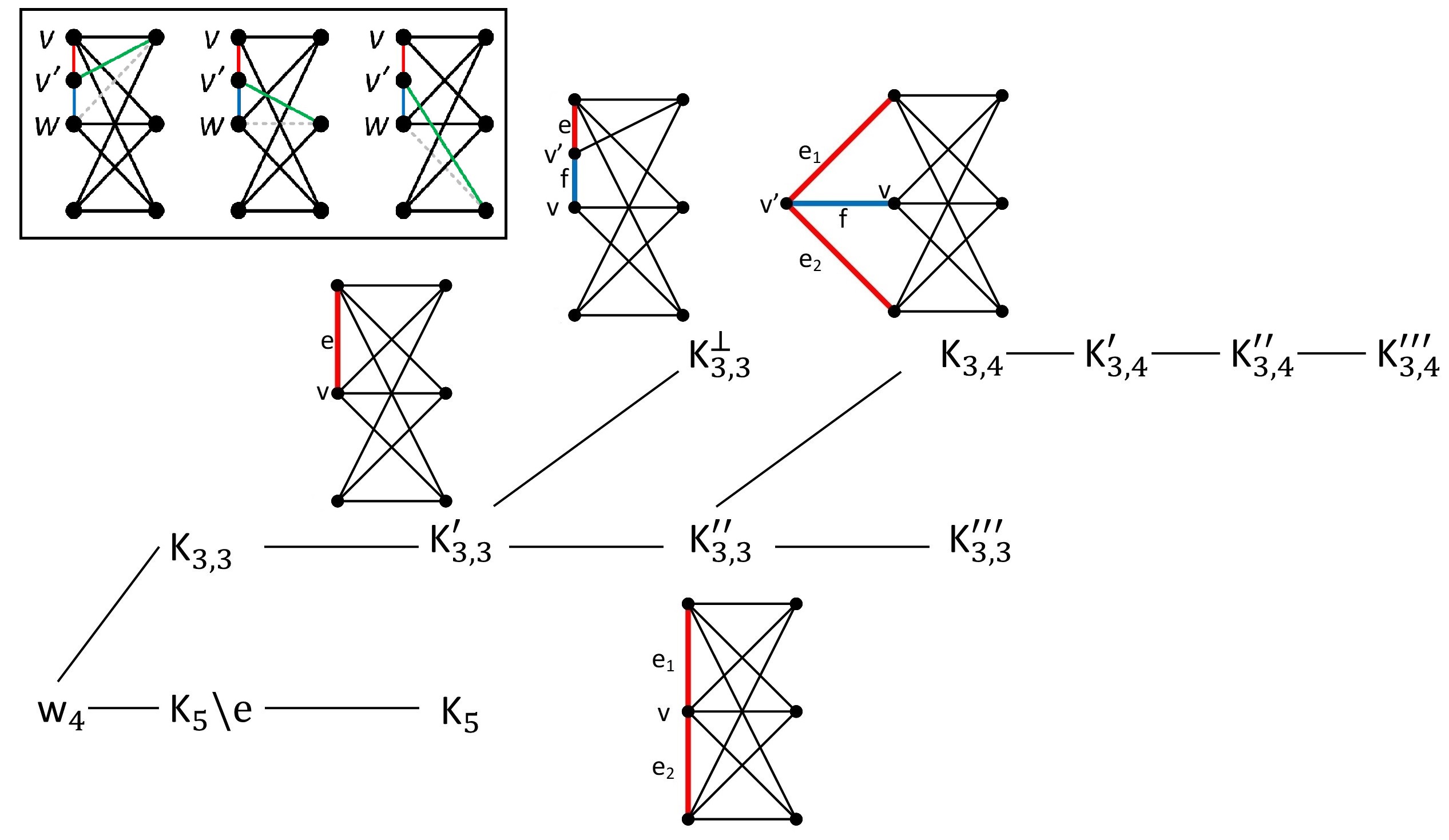}
\caption{ Base case for Dirac's Theorem}
\label{Dirac-basecase}
\end{figure}

\begin{proof} One direction is immediate. The graphs mentioned in the statement of the theorem are in $\mathcal M$ since they do not have two vertex disjoint cycles, and therefore no prism minor. 

Conversely, suppose $G$ is a simple 3-connected graph with no prism minor. Theorem \ref{SplitterTheorem} implies that, except for the wheels, we can construct every 3-connected graph by starting with $W_4$ and performing edge additions and vertex splits (see Figure \ref{W4} in the Appendix). Up to isomorphism, $W_4+e = K_5\backslash e$ and $W_4+\{e_1, e_2\} = K_5$. All graphs on 5 vertices, that is $W_4$, $K_5\backslash e$, and $K_5$ have no prism minor since the prism has 6 vertices. The graph $W_4$ has two non-isomorphic vertex splits, namely $K_{3, 3}$ and the prism graph; $K_5\backslash e$ has two non-isomorphic vertex splits, namely $prism +e$ and $K'_{3, 3}+e$; and $K_5$ has a unique vertex split, namely $K_{3,3}$ with an edge in each vertex class. Since the prism is excluded, we may assume $G$ has a $K_{3,3}$-minor.
The proof is by induction on $n\ge 7$. The base case is illustrated in Figure \ref{Dirac-basecase}. The result holds for $n=7$ since $K_{3,3}^{\perp}$ is the unique vertex split of $K'_{3,3}$ and it has a prism minor. The graph $K^{\perp}_{3,3}$ is obtained by adding edge $e=vw$ to $K'_{3,3}$ and splitting one of its end vertices as specified in Theorem \ref{FCT-BottomUp}(ii). Vertex $v$ is split into $v$ and $v'$, and forming a triad are edges $vv'$, $v'w$, and a third edge (shown in green) that is unlinked from $v$ and linked to $v'$. There are choices for the third edge since any edge in $N(v)-e$ may be unlinked from $v$ and linked to $v'$. However, due to the symmetry of edges in $K_{3,3}$, all choices give graphs isomorphic to $K^{\perp}_{3,3}$. See the three graphs drawn inside the box in Figure \ref{Dirac-basecase}. Thus, for $n=6$, the minimal 3-connected graph with no prism minor is $K_{3,3}$ and the monarch is $K_{3,3}'''$.

Assume that the result hold for 3-connected graphs with no prism minor on $n-1$ vertices, that is, the minimally 3-connected graph on $n-1$ vertices is $K_{3, n-4}$ and the monarch is $K'''_{3, n-4}$. Suppose $G$ has $n$ vertices and first, suppose $G$ is minimally 3-connected. By the induction hypothesis and Theorem \ref{SST-BottomUp}, $G$ can be constructed from $K_{3,n-4}$ using Operations (i), (ii), and (iii). Since $K_{3, n-4}$ is cubic, it has no vertex splits. Using Operation (ii) and the symmetry of the edges present in $K'_{3,n-4}$, the only vertex split of $K'_{3, n-4}$ to consider is $K^{\perp}_{3, n-4}$, and it has a prism minor since $K^{\perp}_{3,3}$ has a prism minor. Operation (iii) gives precisely one graph when it is executed, namely $K''_{3, n-4}+\{e_1, e_2\}\circ f =K_{3, n-3}$. Therefore $G\cong K_{3, n-3}$. 
Second, the only 3-connected edge additions of $K_{3, n-3}$ are $K'_{3, n-3}$, $K''_{3, n-3}$, and $K'''_{3, n-3}$. Observe that $K'''_{3, n-3}+e$ has a minor isomorphic to $K_{3, n-4}+e$, which has a prism minor. Therefore, the monarch is $K_{3, n-3}$. \end{proof}

\begin{figure}[h] 
\centering
\includegraphics[width=6in]{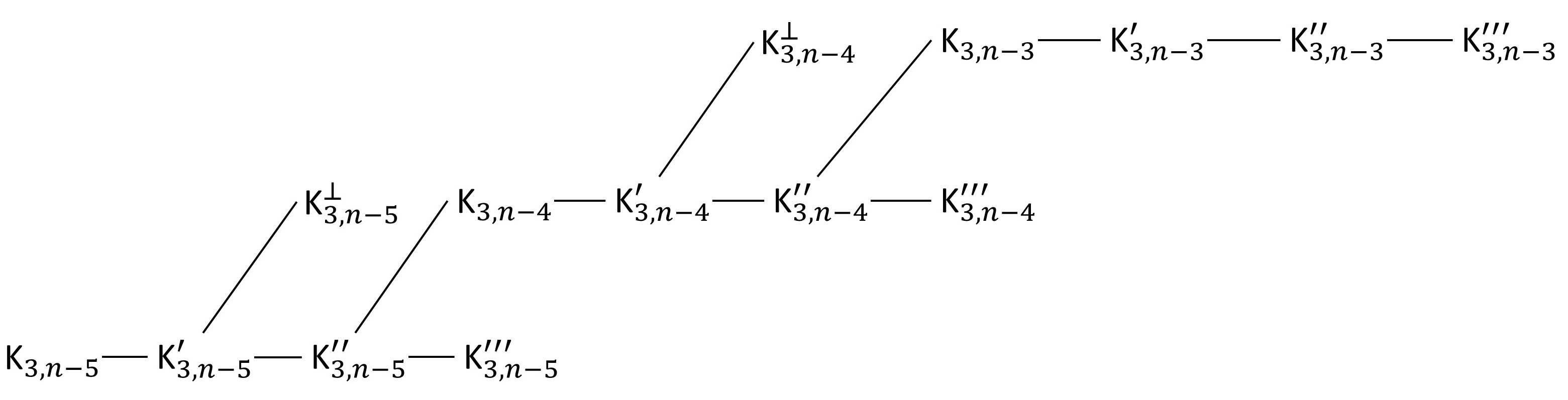}
\caption{Growth pattern of simple 3-connected graphs with no prism minor. } 
\label{Dirac-growthpattern} 
\end{figure}

\subsection{An extension of Halin's theorem } 
$ $

Halin proved that if $G$ is a minimally $3$-connected graph on $n\ge 8$ vertices, then $|E(G)|\le 3n-9$ and equality holds if and only if $G\cong K_{3, n-3}$. We will give a different proof of Halin's result and extend it by proving that if $G$ is minimally 3-connected with $|E(G)|=3n-10$, then $G\cong K^{\perp}_{3, n-4}$. 

Before presenting the proof, let us look at the small cases to understand why the pattern begins at higher values of $n$. Observe that $|E(W_{n-1})|=2(n-1)$. If $n=6$, then $W_5$ has 10 edges, whereas $K_{3, 3}$ and the prism graph have only 9 edges. So the induction argument cannot begin at $n=6$. If $n=7$, then both $W_6$ and $K_{3, 4}$ have 12 edges, and if $n=8$, both $W_7$ and $K^{\perp}_{3, 4}$ have 14 edges. For $n\ge 8$, $3n-9 > 2n-2$ and for $n\ge 9$, $3n-10 > 2n-2$, so the wheels are no longer extremal graphs. We will prove that the infinite families $K_{3, n-3}$ for $n\ge 8$ and $K^{\perp}_{3, n-4}$ for $n\ge 9$ are the unique families of sizes $3n-9$ and $3n-10$, respectively.

\begin{lemma} \label{HalinLemma}Suppose $G$ is a minimally $3$-connected graph with a prism minor, but no $K_{3,3}$-minor, then for $n\ge 9$, $|E(G)|\le 3n-11$.
\end{lemma}

\begin{proof} Theorem \ref{FCT-BottomUp} implies that $G$ may constructed from the prism graph using the three specified operations. As shown in the Appendix (Figure \ref{tightbound-3}) the first time Operation (iii) gives a minimally 3-connected graph is when $n=8$ and the resulting graph is $K^{\perp}_{3, 4}$ with 9 vertices and 14 edges. Clearly, $K^{\perp}_{3, 4}$ has a $K_{3,3}$-minor. So for $n=9$, $|E(G)|\le 13 = 3n-11$. Assume the result holds for $n-1$ vertices and suppose $G$ has $n$ vertices. By Theorem \ref{FCT-BottomUp} and the induction hypothesis, $|E(G)|\le |E(G')|+3 \le 3(n-1)-11+3=3n-11.$ \end{proof}

\begin{theorem} \label{HalinExtension} Let $G$ be a minimally $3$-connected graph on $n\ge 7$ vertices. Then $|E(G)|\le 3n-9$. Moreover, if $n\ge 8$, $|E(G)|= 3n-9$ if and only if $G\cong K_{3, n-3}$ and if $n\ge 9$, $|E(G)|= 3n-10$ if and only if $G\cong K^{\perp}_{3, n-4}$. 
\end{theorem}

\begin{figure}[h]
\centering
\includegraphics[width=5in]{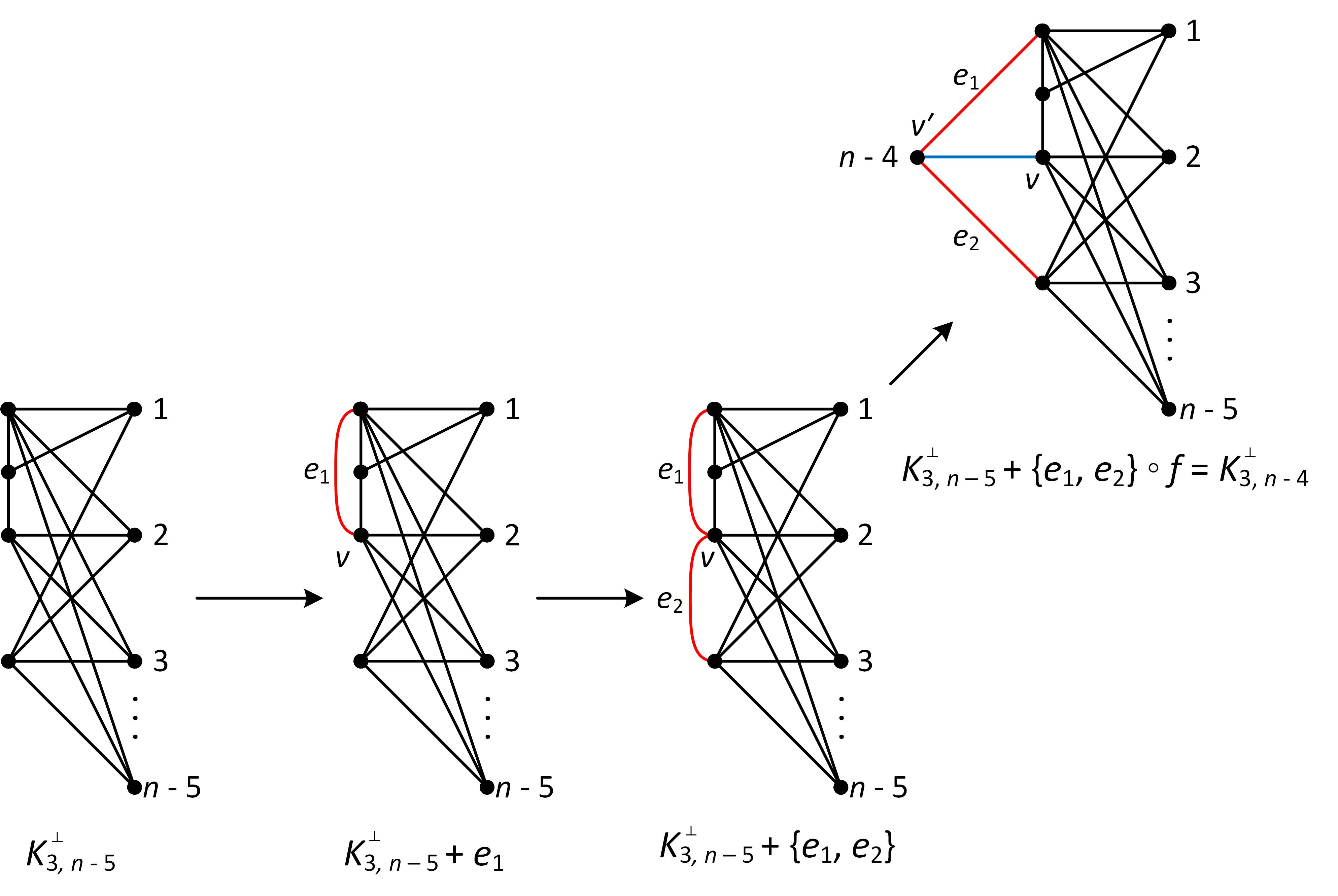}
\caption{ Growth pattern of minimally 3-connected graphs with $3n-10$ edges.}
\label{Halin-growthpattern}
\end{figure}

\begin{proof}Suppose $G$ is a minimally 3-connected graph with $n\ge 7$ vertices. If $G\cong W_{n-1}$, then 
$|E(G)|\le 3n-9$, so we may assume that $G\not\cong W_{n-1}$. Since $K_{3,3}$ and the prism graph are the only vertex splits of $W_4$, Theorem \ref{SplitterTheorem} implies that $G$ has a $K_{3, 3}$-minor or a prism-minor. 
If $G$ has a prism-minor, but no $K_{3,3}$-minor, then $|E(G)|\le 3n-11$ by Lemma \ref{HalinLemma}. Therefore, suppose $G$ has a $K_{3,3}$-minor. Theorem \ref{FCT-Bound} implies that $$|E(G)| \le 9+ 3(|V(G)|-6)=3|V(G)|-9.$$ We will prove by induction on $n\ge 7$ that if $G$ has $3n-9$ edges, then $G\cong K_{3, n-3}$ and if $G$ has $3n-10$ edges, then $G\cong K^{\perp}_{3, n-4}$. Suppose $n=7$. Figure \ref{Dirac-basecase} illustrates how $K_{3,3}+\{e_1, e_2\}\circ f=K_{3,4}$ and Figure \ref{FCT-BottomUp-Example2} illustrates how $K^{\perp}_{3,3}+\{e_1, e_2\}\circ f=K^{\perp}_{3,4}$. In Figure \ref{FCT-BottomUp-Example2}, two edges $e_1$ and $e_2$ are added to a common vertex in $K^{\perp}_{3,3}$. When $v$ is split to obtain $K^{\perp}_{3,3}+\{e_1, e_2\} \circ f$ the first graph leads to $K^{\perp}_{3,4}$ as required. In the second graph, when the common vertex $v$ is split, $K^{\perp}_{3,3}+\{e_1, e_2\} \circ f$ is not minimally 3-connected. This is because in the second graph $v$ is adjacent to $d$ which has degree 4 vertex, so edge $dv$ (shown in green) is deletable in $K^{\perp}_{3,3}+\{e_1, e_2\} \circ f$. In the third graph, the two edges $e_1$ and $e_2$ are incident to the end vertices of the same edge (shown in green) so that edge is deletable in $K^{\perp}_{3,3}+\{e_1, e_2\} \circ f$.

\begin{figure}[h]
\centering
\includegraphics[width=2.5in]{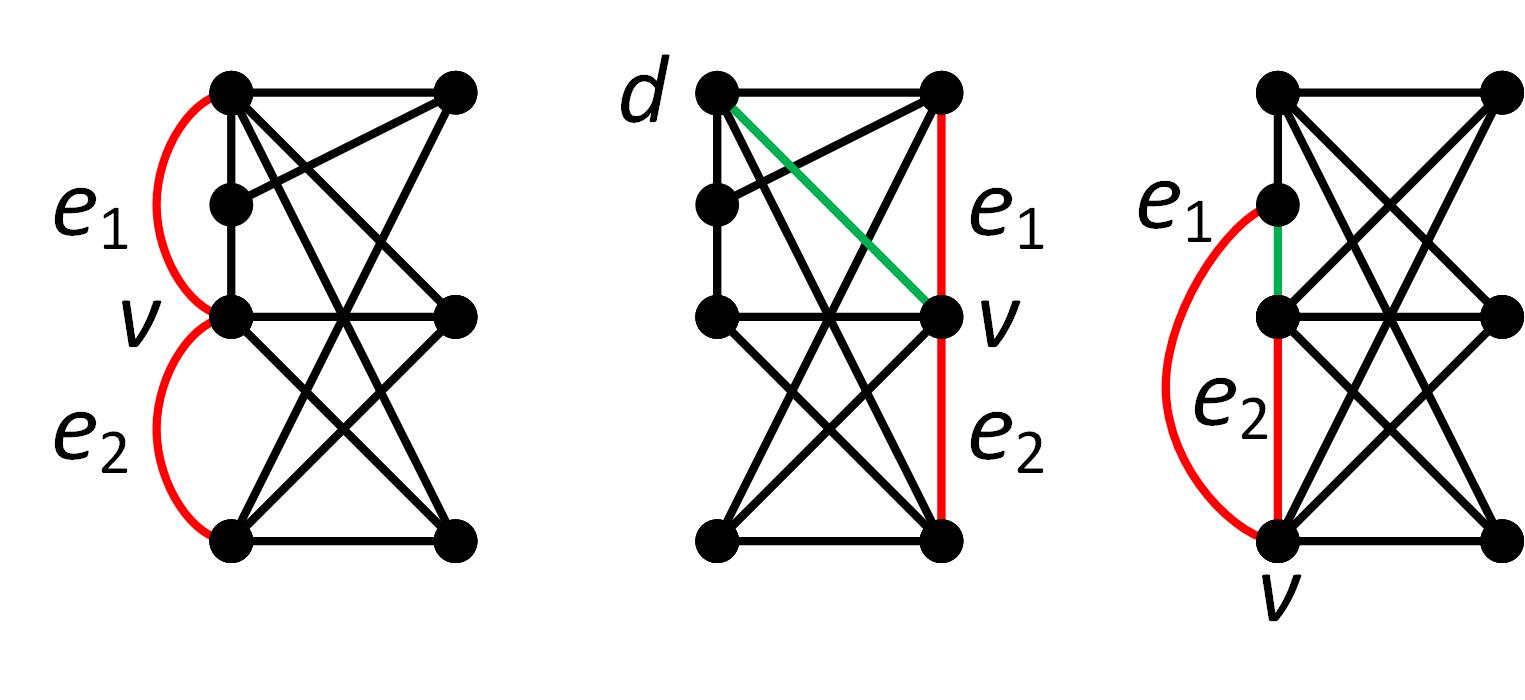}
\caption{ Example showing possibilities for $K^{\perp} _{3,3}+\{e_1, e_2\}$.}
\label{FCT-BottomUp-Example2}
\end{figure}

Assume that the result is true for $n-1$ vertices, that is, $K_{3, n-4}$ and $K^{\perp}_{3, n-5}$ are, respectively, the unique minimally 3-connected families with $3(n-1)-9$ edges and $3(n-1)-10$ edges. 

Suppose $G$ has $n$ vertices and $3n-9$ edges. Since $G$ has one more vertex and three more edges than $K_{3, n-4}$, the induction hypothesis and Theorem \ref{FCT-BottomUp} imply that $G$ is obtained from $K_{3, n-4}$ using Operation (iii). In this case $K_{3, n-4}+\{e_1, e_2\} \circ f=K_{3, n-3}$. 
Suppose $G$ has $n$ vertices and $3n-10$ edges. Here we have two possibilities. 
\begin{itemize}
\item [(i)] $G$ has one more vertex and one more edge than $K'_{3,n-4}$ which has $n-1$ vertices and $3(n-4)+1=3n-11$ edges. Using the same argument as in the proof of Theorem \ref{DiracTheorem1963}, adding an edge $e$ to $K'_{3,n-4}$ and splitting an end vertex of $e$ gives the unique graph $K^{\perp}_{3, n-4}$. That is $K'_{3,n-4}+e\circ f = K^{\perp}_{3, n-4}$. 
\item [(ii)]$G$ has one more vertex and three more edges than $K^{\perp}_{3, n-5}$ which has $n-1$ vertices and $3(n-5)+2=3n-13$ edges. The induction hypothesis and Theorem \ref{FCT-Bound} imply that $G$ may be obtained from $K^{\perp}_{3, n-5}$ using Operation (iii), as shown in Figure \ref{Halin-growthpattern}. In this case $K^{\perp}_{3, n-5}+\{e_1, e_2\} \circ f=K^{\perp}_{3, n-4}$, . 
\end{itemize}
Thus $K^{\perp}_{3, n-4}$ is the unique graph with $3n-10$ edges. \end{proof}


\section{\bf Conclusion and next steps}

The main result in this paper, presented in both top-down format (Theorem \ref{FCT-TopDown}) and constructive format (Theorem \ref{FCT-BottomUp}) has several applications and we anticipate finding more such applications especially since Dirac's Theorem (Theorem \ref{DiracTheorem1963}) is used in several excluded minor results. It may be worth examining all those results to see how they can also be shortened. Original proofs are often long and complicated and seeking simple proofs is a worthy goal in and of itself. Sometimes, as in the case of Theorem \ref{HalinExtension}, new results hidden in plain sight emerge when we look at them the right way.

We conjecture that planar minimally 3-connected graphs (except for $W_{n-1}$) can be constructed from the prism graph using only Operations (i) and (ii) in Theorem \ref{FCT-BottomUp}. Every candidate for a counterexample to this conjecture turned out to be constructed from some minimally 3-connected minor using only the two operations. If this conjecture is true, then Ota's theorem in \cite{Ota2004} that the maximum number of edges in an $n$-vertex planar minimally 3-connected graph is $2n-3$ would follow immediately.

Often an excluded minor class of graphs is too large to identify all its 3-connected members. In this case the connectivity is raised to 4-connectivity or cyclic 4-connectivity. We will explore extending Theorem \ref{FCT-TopDown} to higher levels of connectivity in future work.


\section*{\bf Appendix} 

\begin{figure}[h]
\centering
\includegraphics[width=4in]{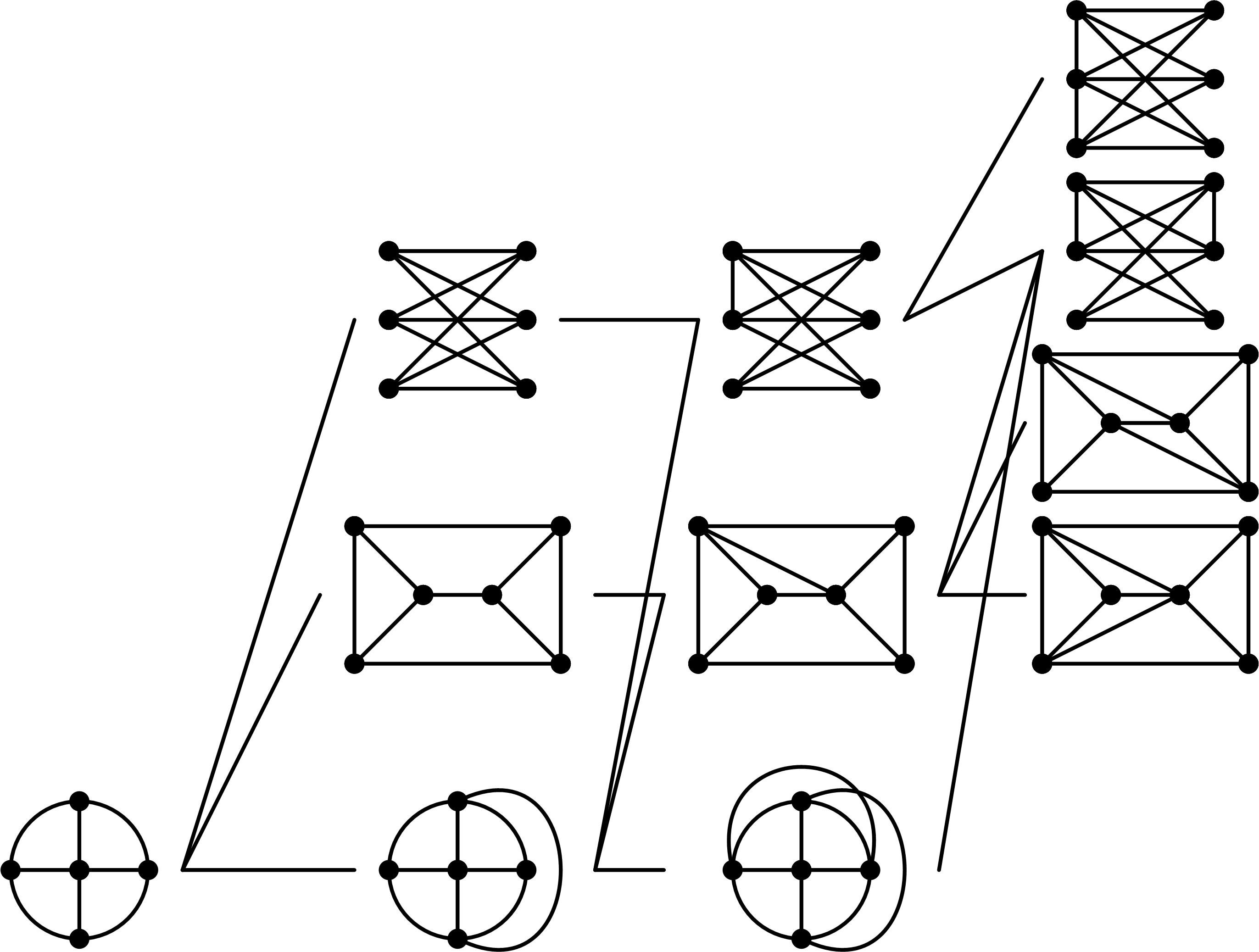}
\caption{Small 3-connected graphs
\label{W4}}
\end{figure}

The appendix has the computational details omitted in the proofs. Figure \ref{W4} illustrates the edge additions and vertex splits of the 3-connected graphs with $n\le 6$ vertices and $m\le 11$ edges, starting with $W_4$. A horizontal line indicates the graph is an edge addition and a slanted line indicates the graph is a vertex split.

\begin{figure}[h]
\centering
\includegraphics[width=5.7in]{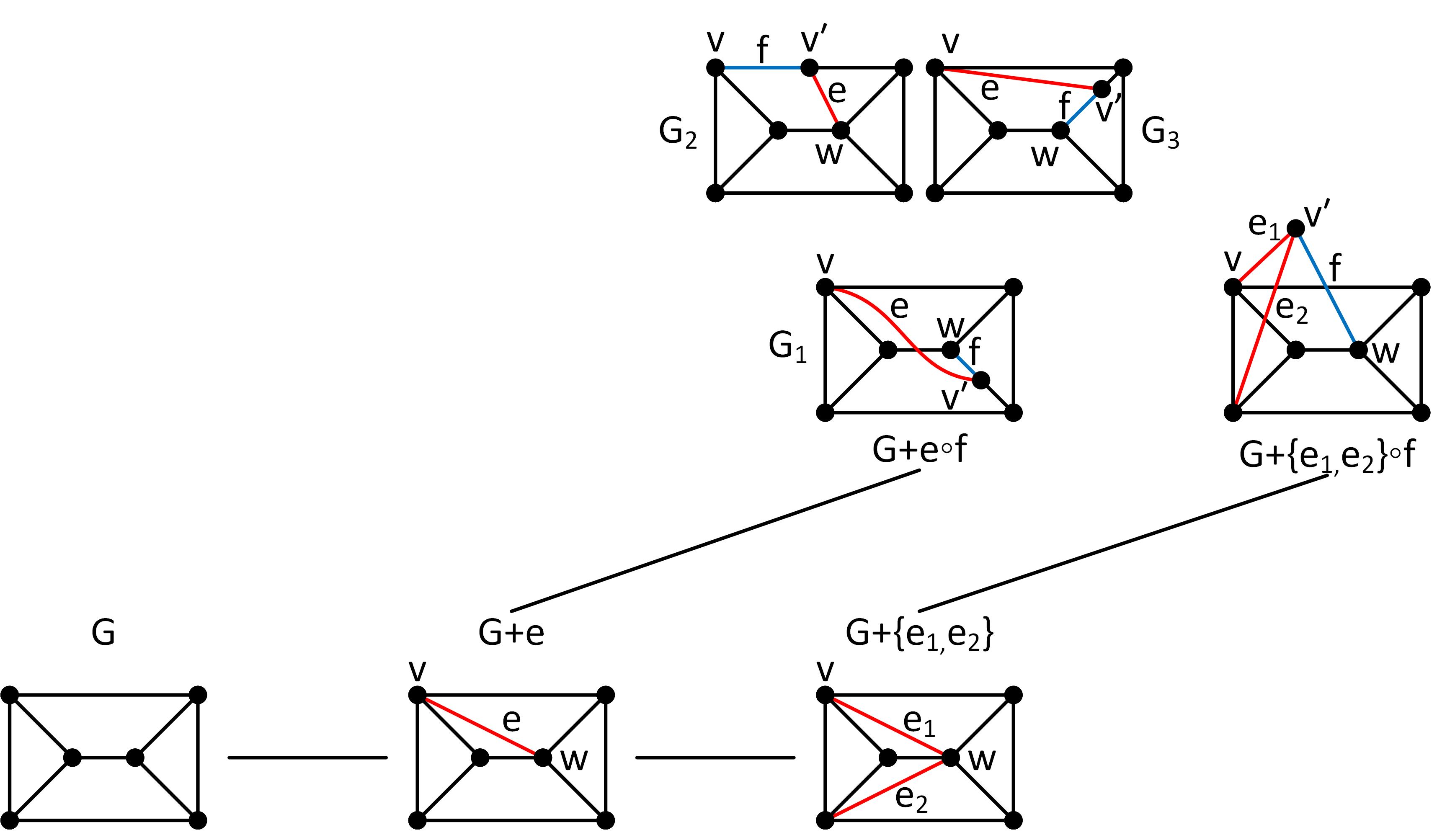}
\caption{The three operations in Theorem \ref{FCT-BottomUp} applied to the prism graph
\label{prismFCT}}
\end{figure}

Figure \ref{prismFCT} illustrates the operations in Theorem \ref{FCT-BottomUp} when $H$ is the prism graph. 
The prism is a cubic graph and therefore has no vertex splits. Up to isomorphism, the prism graph has just one edge addition and although it has three double edge additions, there is only one where both edges are incident to the same vertex. Up to isomorphism, $prism + e_1$ has three vertex splits satisfying the second condition. The single graph that satisfies the third condition is not minimally 3-connected.

\begin{figure}[h]
\centering
\includegraphics[width=4in]{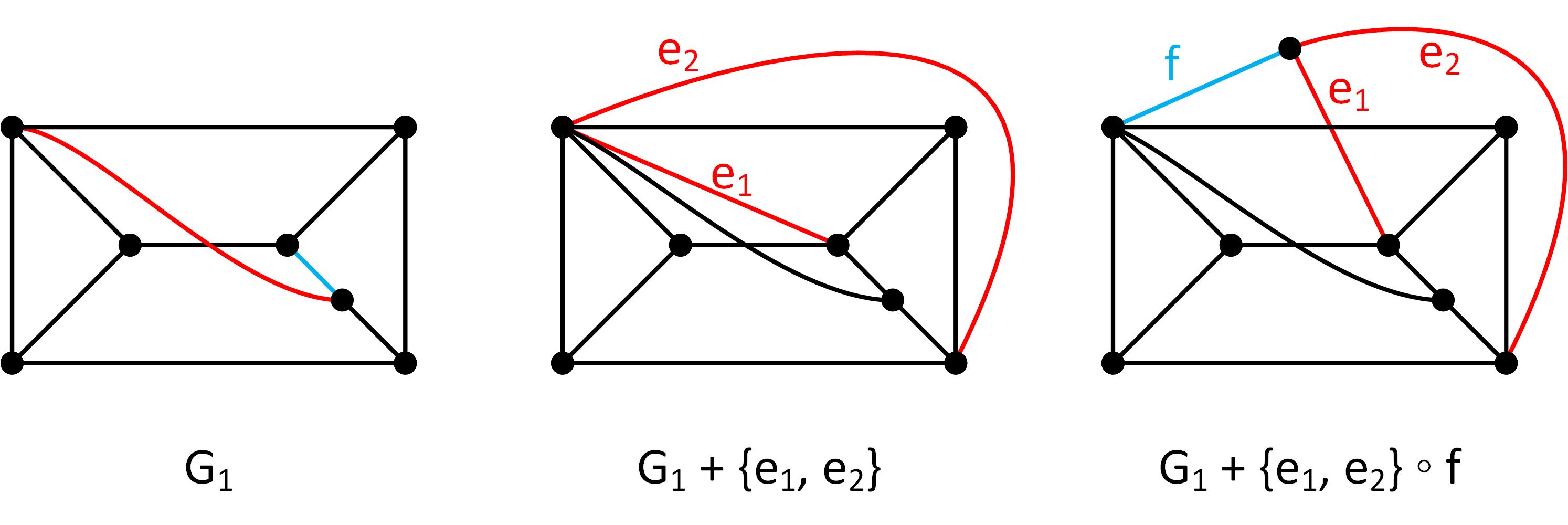}
\caption{First time Theorem \ref{FCT-BottomUp}(iii) is used when $H$ is the prism graph}
\label{tightbound-3}
\end{figure}

Construction of all minimally 3-connected graphs with 8 vertices and 14 edges starting with the prism graph shows that the first time that Operation (iii) is used is when $n=8$ vertices. Applying Theorem \ref{FCT-BottomUp}(iii) to $G_2$ and $G_3$ gives a graph that is not minimally 3-connected. However, applying Theorem \ref{FCT-BottomUp}(iii) to $G_1$ gives the graph shown Figure \ref{tightbound-3}. Observe that $G_1+\{e_1, e_2\}\circ f \cong K^{\perp}_{3, 4}$, so it is a non-planar graph. 

$ $


\end {document}